\documentclass[preprint]{imsart}

\usepackage{amsthm,amsmath}

\usepackage{amssymb}
\usepackage{dsfont}
\usepackage{amsfonts}

\RequirePackage[colorlinks,citecolor=blue,urlcolor=blue,breaklinks]{hyperref}
\RequirePackage[authoryear]{natbib}

\startlocaldefs
\newtheorem{theorem}{Theorem}

\newtheorem{hyp}[theorem]{Hypothesis}

\theoremstyle{remark}
\newtheorem{remark}{Remark}
\newtheorem{ex}[remark]{Exemple}

\endlocaldefs

\newcommand{\prob}[1]{\mathbb{P}\big{\{}#1\big{\}}}
\newcommand{\esp}[1]{\mathbb{E}\big{[}#1\big{]}}
\newcommand{\ind}[1]{\mathbf{1}_{\{ #1 \}}}

\def\R{\mathbb{R}}
\def\B{\mathbb{B}}
\def\N{\mathbb{N}}

\def\E{\mathbb{E}}

\begin{document}

\begin{frontmatter}

\title{Adiabatic reduction of models of stochastic gene expression with bursting}
\runtitle{adiabatic reduction for PDMP}

\begin{aug}
\author{\fnms{Romain} \snm{Yvinec}\ead[label=e1]{yvinec@math.univ-lyon1.fr}\corref{}\thanksref{t1}}
\thankstext{t1}{corresponding author}
\runauthor{R. Yvinec1}
\address{Universit\'e de Lyon CNRS UMR 5208 Universit\'e Lyon 1\\ Institut Camille Jordan
43 blvd. du 11 novembre 1918 \\ F-69622 Villeurbanne Cedex France \\ \printead{e1}}
\end{aug}

\begin{abstract}
This paper considers adiabatic reduction in both discrete and continuous models of stochastic gene expression. In gene expression model, the concept of bursting is a production of several molecules simultaneously and is generally represented as a jump terms of random size. In a general two-dimensional birth and death discrete model, we prove that under specific assumptions and scaling (that are characteristics of the mRNA-protein system) an adiabatic reduction leads to a one-dimensional discrete-state space model with bursting production. The burst term appears then through the reduction of the first variable. In a two-dimensional continuous model, we also prove that an adiabatic reduction can be performed in a stochastic slow/fast system. In this gene expression model, the production of mRNA (the fast variable) is assumed to be bursty and the production of protein (the slow variable) is linear as a function of mRNA. When the dynamics of mRNA is assumed to be faster than the protein dynamics (due to a mRNA degradation rate larger than for the protein) we prove that, with the appropriate scaling, the bursting phenomena can be transmitted to the slow variable. We show that  the reduced equation is either a stochastic differential equation with a jump Markov process or a deterministic ordinary differential equation depending on the scaling that is appropriate. 

These results are significant  because adiabatic reduction techniques seem to have not been applied to a stochastic differential system containing  a jump Markov process. Last but not least, for our particular system, the adiabatic reduction allows us to understand what are the necessary conditions for the bursting production-like of protein to occur.
\end{abstract}



\end{frontmatter}

\section*{Introduction}

The adiabatic reduction techniques give results that allow to reduce the dimension of a system and justify the use of an effective set of reduced equations in lieu of dealing with a full, higher dimensional model, if different time scales occur in the system. Adiabatic reduction results for deterministic systems of ordinary differential equations have been available since the very precise results of \cite{Tikhonov1952} and \cite{Fenichel1979}. The simplest results, in the hyperbolic case, give an effective construction of an uniformly asymptotically  stable slow manifold (and hence a reduced equation) and prove the existence of an invariant manifold near the slow manifold, with (theoretically) any order of approximation of this invariant manifold. Such precise and geometric results have been generalized to random systems of stochastic differential equation with Gaussian white noise (\cite{Berglund2006}, see also \cite{Gardiner1985} for previous work on the Fokker-Planck equation). However, to the best of our knowledge, analogous  results for stochastic differential equations with jumps have not been obtained.

The present paper gives a theoretical justification of an adiabatic reduction of a particular piecewise deterministic Markov process \citep{Davis1984}. The results we obtain do not give a bound on the error of the reduced system, but they do allow us to justify the use of a reduced system in the case of a piecewise deterministic Markov process. In fact, we prove limit theorems using martingale strategy, in a similar manner than in recent papers such as \cite{Radulescu}, \cite{Kang} and \cite{Riedler2012}, where general convergence results for discrete models of stochastic reaction networks are given. In particular, these papers give alternative scaling of the traditional ordinary differential equation and the diffusion approximation depending on the different scaling chosen (see \cite{Ball2006} for some examples in a reaction network model). After the scaling, the limiting models can be deterministic (ordinary differential equation), stochastic (jump Markov process), or hybrid (piecewise deterministic process). For illustrative and motivating examples given by a simulation algorithm, see \cite{Haseltine2002,Rao2003,Goutsias2005}. However, we emphasize that we do not consider here a continuous approximation of a discrete model. Rather, we perform adiabatic reduction on both discrete state-space and continuous state-space  models. Time-scale reduction have been considered in \cite{Kang}, but not on the kind we perform here.

Our particular model is meant to describe stochastic gene expression with explicit bursting \citep{Friedman2006}. In discrete state-space bursting models, the variables evolve under the action of a discrete birth and death process, interrupted by discrete positive jumps of random sizes. In continuous state-space bursting models,  the variables evolve under the action of a continuous deterministic dynamical system, interrupted by positive jumps of random sizes. In both cases, the positive jumps model the burst production of \textit{several} molecules instantaneously. In that sense, the convergence theorems we obtain in this paper can be seen as an example in which there is a reaction with size between $0$ and $\infty$. We hope that the results here are generalizable to give insight into adiabatic reduction methods in more general stochastic hybrid systems \citep{Hespanha2006,Bujorianu2004}. We note also that more geometrical approaches have been proposed to reduce the dimension of such systems in \cite{Bujorianu2008}.

Biologically, the bursting of mRNA or protein molecules is defined as the production of several molecules within a very short time, \textit{indistinguishable} within  the time scale of the measurement.  In the biological context of models of stochastic gene expression, explicit models of bursting mRNA and/or protein production have been analyzed recently, either using a discrete  \citep{Shahrezaei2008,Lei2009} or a continuous formalism \citep{Friedman2006,Mackey2011} as more and more experimental evidence from single-molecule visualization techniques has revealed the ubiquitous nature of this phenomenon  \citep{Ozbudak2002,Golding2005,Raj2006,Elf2007,Xie2008,Raj2009,Suter2011}. Traditional models of gene expression are composed of \textit{at least} two variables (mRNA and protein, and sometimes the DNA state). The use of a reduced one-dimensional model (that has the advantage that it can be solved analytically) has been justified so far by an argument concerning the stationary distribution in \cite{Shahrezaei2008}. However, it is clear that two different models may have the same stationary distribution but very different behavior (continuous or discontinuous trajectories, monostable or bistable, etc; for an example in that context, see \cite{Mackey2011}). Hence, our results are of importance to rigorously prove the validity of using a reduced model. Our results are based on the standard assumption that the mRNA molecules have a shorter lifetime than the protein molecules, that is widely observed in both prokaryotes and eukaryotes (\cite{Schwanhausser2011}). Depending on the assumed scaling of other kinetic parameters within the mRNA degradation rates, different limiting models are obtained.

The paper is organized as follows. In the first section, we prove a reduction results for a discrete state-space model, that is a two-dimensional birth and death process. Assumptions on the birth and death rates are in agreement with a standard model of gene expression for the mRNA-protein system. That is both variables remain positive and birth of the second variable can occur only if the first variable is positive. Using an appropriate scaling of birth and death rates, we prove that this model converges to a general one-dimensional discrete bursting model.

In the second section, we prove a reduction for a continuous state-space model, that is a two-dimensional piecewise deterministic model of gene expression with a jump production term for the first variable. Using appropriate scaling on parameters, we prove that this model converge either to a deterministic ordinary differential equation or to a one-dimensional continuous bursting model.

\section{A bursting model from a two-dimensional discrete model}\label{ssec:reduction_2d}

The fact that bursting models arise as a reduction procedure of a higher dimensional model was already observed in \cite{Shahrezaei2008} and \cite{Radulescu}. In \cite{Shahrezaei2008}, the authors show that, within an appropriate scaling, the stationary distribution of a 2-dimensional discrete model converge to the stationary distribution of a 1-dimensional bursting model. The authors used analytic methods through the transport equation on the generating function. Their result seems to be restricted to first-order kinetics. The first variable is a fast variable that induces infrequent kicks to the second one.
In \cite{Radulescu}, the authors show that, within an appropriate scaling, a fairly general discrete state space model with a binary variable converge to a bursting model with continuous state space. The authors obtained a convergence in law of the solution through Martingale techniques. The binary variable is a fast variable that, when switching in an "ON" state, induces kicks to the other variable.

We present below analogous result of \cite{Radulescu} when the fast variable is similar to the one of \cite{Shahrezaei2008}. Our limiting model is still a discrete state space model. These results are more precise than the one of \cite{Shahrezaei2008}, and more general (some kinetics rates can be non-linear). We use martingales techniques, with a proof that is similar to \cite{Radulescu} and also inspired by results from \cite{Kang}.  We present below the model, then state our result in the subsection \ref{ssec:result_disc}, and divide the proof in the three next subsections  \ref{ssec:proof_disc1}-\ref{ssec:proof_disc3}.

We consider the following two-dimensional stochastic kinetic chemical reaction model
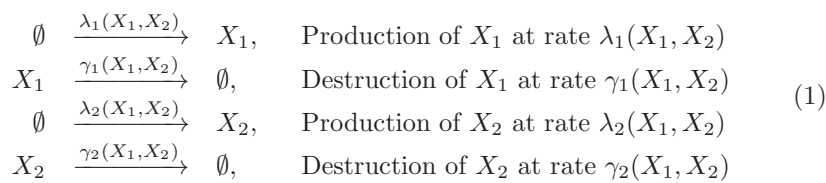
\begin{equation}\label{modelD2}
\begin{array}{rcll}

\emptyset  &\xrightarrow{\lambda_1(X_1,X_2)}&  X_1, & \,\,\,\,\text{Production of }X_1 \text{ at rate }\lambda_1(X_1,X_2)
\\
X_1  &\xrightarrow{\gamma_1(X_1,X_2)}&  \emptyset, & \,\,\,\,\text{Destruction of }X_1\text{ at rate }\gamma_1(X_1,X_2)
\\
\emptyset  &\xrightarrow{\lambda_2(X_1,X_2)}&  X_2, &\,\, \,\,\text{Production of }X_2\text{ at rate }\lambda_2(X_1,X_2)
\\
X_2  &\xrightarrow{\gamma_2(X_1,X_2)}& \emptyset,  & \,\,\,\,\text{Destruction of }X_2\text{ at rate }\gamma_2(X_1,X_2)
\end{array}
\end{equation}

with $\gamma_1(0,X_2)=\gamma_2(X_1,0)=0$ to ensure positivity. This model can be represented by a continuous time Markov chain in $\N^2$, and is then a general birth and death process in $\N^2$. It can be described by the following set of stochastic differential equations 
 \begin{align*}
X_1(t)&=X_1(0)+Y_1\Big{(}\int_0^t\lambda_1(X_1(s),X_2(s))ds\Big{)}-Y_2\Big{(}\int_0^t\gamma_1(X_1(s),X_2(s))ds \Big{)},\\
X_2(t)&=X_2(0)+Y_3\Big{(}\int_0^t\lambda_2(X_1(s),X_2(s))ds\Big{)}-Y_4\Big{(}\int_0^t\gamma_2(X_1(s),X_2(s))ds \Big{)},
 \end{align*}
where $Y_i$, for $i=1...4$ are independent standard poisson processes. The generator of this process is given by
\begin{equation}
\begin{aligned}
\mathbb{B} f(X_1,X_2)= & \lambda_1(X_1,X_2)\Big{[}f(X_1+1,X_2)-f(X_1,X_2)\Big{]}\\
& +\gamma_1(X_1,X_2)\Big{[}f(X_1-1,X_2)-f(X_1,X_2)\Big{]} \\
& +\lambda_2(X_1,X_2)\Big{[}f(X_1,X_2+1)-f(X_1,X_2)\Big{]} \\
& +\gamma_2(X_1,X_2)\Big{[}f(X_1,X_2-1)-f(X_1,X_2)\Big{]},
\end{aligned}
\end{equation}
for every bounded function $f$ on $\N^2$. 

\begin{ex}\label{ex:1}
We have in mind the standard mRNA-Protein system given by the following choice: $\gamma_i(X_1,X_2)=g_i X_i$ with $g_i>0$  for $i=1,2$, $\lambda_1(X_1,X_2)=\lambda_1(X_2)$ and $\lambda_2(X_1,X_2)=k_2 X_1$ with $k_2>0$. Note however that even in the context of models of gene expression, different models have been proposed, that includes nonlinear feedback of mRNA and/or nonlinear degradation terms \cite{Bose}.
\end{ex}

\subsection{Statement of the result}\label{ssec:result_disc}

We suppose the following scaling holds
\begin{align*}
 \gamma_1^N(X_1,X_2) &= N\gamma_1(X_1,X_2) \\
\lambda_2^N (X_1,X_2) &= N\lambda_2 (X_1,X_2)
\end{align*}
where $N \to \infty$ that is degradation of $X_1$ and production of $X_2$ occurs at a faster time scale than the two other reactions. Then $X_1$ is degraded very fast, and induces also as a very fast production of $X_2$. The rescaled model is given by

\begin{equation}\label{eq:sderescale}
\begin{aligned}
X_1^N(t)&=X_1^N(0)+Y_1\Big{(}\int_0^t\lambda_1(X_1^N(s),X_2^N(s))ds\Big{)}-Y_2\Big{(}\int_0^tN\gamma_1(X_1^N(s),X_2^N(s))ds \Big{)},\\
X_2^N(t)&=X_2^N(0)+Y_3\Big{(}\int_0^tN\lambda_2(X_1^N(s),X_2^N(s))ds\Big{)}-Y_4\Big{(}\int_0^t\gamma_2(X_1^N(s),X_2^N(s))ds \Big{)},
 \end{aligned}
\end{equation}
and the generator of this process is given by
\begin{equation}\label{eq:generator_scale}
\begin{aligned}
  \mathbb{B}_N f(X_1,X_2)= & \lambda_1(X_1,X_2)\Big{[}f(X_1+1,X_2)-f(X_1,X_2)\Big{]}\\
  &+N\gamma_1(X_1,X_2)\Big{[}f(X_1-1,X_2)-f(X_1,X_2)\Big{]} \\
&+N\lambda_2(X_1,X_2)\Big{[}f(X_1,X_2+1)-f(X_1,X_2)\Big{]}\\
&+\gamma_2(X_1,X_2)\Big{[}f(X_1,X_2-1)-f(X_1,X_2)\Big{]}.
\end{aligned}
\end{equation}

We can prove the following reduction holds:
\begin{theorem}\label{thm:adiabD2}
We assume that 
\begin{enumerate}
  \item The degradation function on $X_2$ satisfies $\gamma_2(X_1,0)\equiv 0$.
  \item The degradation function on $X_1$ satisfies $\gamma_1(0,X_2)\equiv 0$, and
  \begin{equation*}
    \inf_{X_1\geq 1, X_2\geq 0} \gamma_1(X_1,X_2)=\underline{\gamma}>0.
  \end{equation*}
  \item The production rate of $X_2$ satisfies $\lambda_2(0,X_2)=0$.
  \item The production rate function $\lambda_1$ and $\lambda_2$ are linearly bounded by $X_1+X_2$.
  \item Either $\lambda_1$ or $\lambda_2$ is bounded.
\end{enumerate}
Let $(X_1^N,X_2^N)$ the stochastic process whose generator is $\B_N$ (defined in eq.~\eqref{eq:generator_scale}). Assume that the initial vector $(X_1^N(0),X_2^N(0))$ converges in distribution to $(0,X(0))$, as $N\to\infty$. Then, for all $T>0$, $(X_1^N(t),X_2^N(t))_{t \geq 0}$ converges in $L^1(0,T)$ (and in $L^p$, $1 \leq p <\infty$) to $(0,X(t))$ where $X(t)$ is the stochastic process whose generator is given by
\begin{multline}\label{eq:generator_limit}
 \B_{\infty}\varphi(X)=\lambda_1(0,X)\Big{(}\int_0^{\infty}P_t(\gamma_1(1,.)\varphi(.))(X)dt-\varphi(X)\Big{)}\\ +\gamma_2(0,X)\Big{[}\varphi(X-1)-\varphi(X)\Big{]} ,
\end{multline}
where
\begin{equation*}
 P_tg(X)=\esp{g(Y(t,X)e^{-\int_0^t \gamma_1(1,Y(s,X))ds}},
\end{equation*}
and $Y(t,X)$ is the stochastic process starting at $X$ at $t=0$ whose generator is given by
\begin{equation*}
 Ag(Y)=\lambda_2(1,Y)\big{(}g(Y+1)-g(Y)\big{)}.
\end{equation*}
\end{theorem}
\begin{remark}
The first three hypotheses of theorem \ref{thm:adiabD2} are the main characteristics of the mRNA-protein system (see example~\ref{ex:1}). Basically, they impose that quantities remains non-negative, that the first variable has always the possibility to decrease to $0$ (no matter the value of the second variable), and that the second variable cannot increase when the first variable is $0$. Hence these three hypotheses will guarantee that (with our particular scaling) the first variable converges to $0$, and will lead to an intermittent production of the second variable. The last two hypotheses are more technical, and guarantee that the Markov chain is not explosive, and hence well defined for all $t \geq 0$, and that the limiting model is well defined too.
\end{remark}

\begin{remark}
 The above expression eq.~\eqref{eq:generator_limit} is a generator of a bursting model for a ``general bursting size distribution``. 
 For instance, for linear function $\gamma_1(X_1,X_2)= g_1X_1$, and $\lambda_2(X_1,X_2)=k_2X_1$, we have
\begin{eqnarray*}
 P_t(\gamma_1(.)\varphi(.))(p)&=&g_1P_t(\varphi)(p),\\
&=& g_1\E\Big{[}\varphi(Y_t^y)e^{-g_1 t}\Big{]}, \\
&=& g_1e^{-g_1 t}\sum_{z\geq y}\varphi(z)\prob{Y_t^y=z}, \\
&=& g_1e^{-g_1 t}\sum_{z\geq y}\varphi(z)\frac{(k_2 t)^{z-y}e^{-k_2 t}}{(z-y)!} .\\
\end{eqnarray*}
It follows by integration integration by parts that
\begin{equation*}
 \int_0^{\infty}P_t(\gamma_1(.)\varphi(.))(y)dt=\frac{g_1}{g_1+k_2}\sum_{z \geq 0}\varphi(z+y)\Big{(}\frac{k_2}{k_2+g_1}\Big{)}^z,
\end{equation*}
which gives then an additive geometric burst size distribution of parameter $p=\frac{k_2}{k_2+g_1}$, as expected \cite{Shahrezaei2008}.
\end{remark}

 We divide the proof in three steps: moment estimates, tightness and identification of the limit. 

\subsection{Moment estimates}\label{ssec:proof_disc1}

Because production rates are linearly bounded, it is straightforward that with $f(X_1,X_2)=X_1+X_2$ in eq.~\eqref{eq:generator_scale}, there is a constant $C_N$ (that depends on $N$ and other parameters) such that
\begin{equation*}
  \mathbb{B}_N f(X_1,X_2) \leq C_N (X_1+X_2) .
\end{equation*}
Then $\esp{X_1^N(t)+X_2^N(t)}$ is bounded on any time interval $[0,T]$ and 
$$f(X_1^N(t),X_2^N(t))-f(X_1^N(0),X_2^N(0))-\int_0^t \mathbb{B}_N f(X_1^N(s),X_2^N(s))ds$$
is a $L^1$-martingale.

\subsection{Tightness}\label{ssec:proof_disc2}

Clearly, from the stochastic differential equation on $X_1^N$, we must have \\ $X_1^N(t)\to 0$. We can show in fact that the Lebesgue measure of the set \\ $\{t\leq T: X_1^N(t) \>0 \}$ converges to 0. Indeed, taking $f(X_1,X_2)=X_1$ in eq.~\eqref{eq:generator_scale}, we have
\begin{equation}\label{eq:martX1}
 X_1^N(t)-X_1^N(0)-\int_0^t(\lambda_1(X_1^N(s),X_2^N(s))-N\gamma_1(X_1^N(s),X_2^N(s)))ds 
\end{equation}
is a martingale. Thanks to the lower bound assumption on $\gamma_1$, we have
\begin{equation*}
 \underline{\gamma}\esp{\int_0^t \ind{X_1^N(s)\geq 1}ds}\leq \mathbb{E}\int_0^t\gamma_1(X_1^N(s),X_2^N(s))ds.
\end{equation*}
Then, by the martingale property, we deduce from \eqref{eq:martX1}
\begin{equation}\label{eq:boundX1}
  \underline{\gamma}N\esp{\int_0^t \ind{X_1^N(s)\geq 1}ds}\leq \esp{X_1^N(0)}+\int_0^t\esp{\lambda_1(X_1^N(s),X_2^N(s))}ds.
\end{equation}
Now for $X_2^N$ we obtain from eq.~\eqref{eq:sderescale},
\begin{equation*}
 X_2^N(t) \leq  X_2^N(0)+Y_3\Big{(}\int_0^t N\ind{X_1^N(s)\geq 1}\lambda_2(X_1^N(s),X_2^N(s))ds\Big{)} .
\end{equation*}
Let us now distinguish between the two cases. 

\begin{itemize}
 \item Suppose first that $\lambda_2$ is bounded (say by $K$). Then
\begin{equation*}
 \esp{X_2^N(t)} \leq  \esp{X_2^N(0)}+K N\esp{\int_0^t\ind{X_1^N(s)\geq 1}ds}.
\end{equation*}
As $\lambda_1$ is linearly bounded (say by $K$) by $X_1^N+X_2^N$, the upper bound eq.~\eqref{eq:boundX1} becomes
\begin{equation*}
  \underline{\gamma}N\esp{\int_0^t \ind{X_1^N(s)\geq 1}ds}\leq \esp{X_1^N(0)}+K\int_0^t\Big{(}\esp{X_1^N(s)}+\esp{X_2^N(s)}\Big{)}ds.
\end{equation*}
Finally, with eq.~\eqref{eq:martX1}, it is clear that
\begin{equation*}
 \esp{X_1^N(t)} \leq  \esp{X_1^N(0)}+K\int_0^t\Big{(}\esp{X_1^N(s)}+\esp{X_2^N(s)}\Big{)}ds.
\end{equation*}
Hence, with the three last inequalities, we can conclude by the Gr\"onwall lemma that $\esp{X_2^N(t)}$ is bounded on $[0,T]$, \textit{uniformly} in $N$. Then
$$ N\esp{\int_0^T \ind{X_1^N(s)\geq 1}ds}$$
is bounded and $X_1^N\to 0$ in $L^1([0,T],\N)$. By the law of large number, $\frac{1}{N}Y_3(N)$ is almost surely convergent, and hence almost surely bounded. We deduce then there exists a random variable $C$ such that
\begin{equation*}
 X_2^N(t) \leq  X_2^N(0)+N C\int_0^t \ind{X_1^N(s)\geq 1}ds,
\end{equation*}
almost everywhere. By Gr\"{o}nwall lemma and Markov inequality
\begin{equation*}
 \prob{\sup_{t \in [0,T]} X_2^N(t) \geq M } \to 0
\end{equation*}
as $ M \to \infty$, uniformly in $N$. 

\item Now suppose $\lambda_1$ is bounded (say $K$). By the martingale eq.~\eqref{eq:martX1} (and the same lower bound hypothesis on $\gamma_1$, it is clear that 
$$ N\esp{\int_0^T \ind{X_1^N(s)\geq 1}ds}$$
is bounded and $X_1^N\to 0$ in $L^1([0,T],\N)$. Now, let us denote $U^N(t)=\frac{1}{N}X_1^N(t)$, $V^N=\frac{1}{N}X_2^N(t)$ and $W^N= N\ind{X_1^N(t)\geq 1}$ (which is then bounded in $L^1([0,T[)$). From eq.~\eqref{eq:sderescale}, and from the linear bound on $\lambda_2$ (say by $K$)
\begin{equation*}
 V^N(t) \leq  V^N(0)+ \frac{1}{N}Y_3\Big{(}\int_0^t N K W^N (U^N(s)+V^N(s))ds\Big{)} .
\end{equation*}
Then, still by the law of the large number there exists a random variable $C$ such that
\begin{equation*}
 V^N(t) \leq  V^N(0)+ C\int_0^t W^N (U^N(s)+V^N(s))ds,
\end{equation*}
and hence
\begin{equation*}
 X_2^N(t) \leq  X_2^N(0)+ C\int_0^t W^N (X_1^N(s)+X_2^N(s))ds.
\end{equation*}
By Gr\"{o}nwall lemma,
\begin{equation*}
 \sup_{[0,T]} X_2^N(t) \leq  (X_1^N(0)+X_2^N(0))\exp \Big{(}C\int_0^t W^N(s)ds\Big{)},
\end{equation*}
which is then bounded, uniformly in $N$.
\end{itemize}
For any subdivision of $[0,T]$, $0=t_0 <t_1<\cdots <t_n=T$,
\begin{eqnarray*}
 \sum_{i=0}^{n-1} \mid X_2^N(t_{i+1})-X_2^N(t_i) \mid &\leq & \sum_{i=0}^{n-1} Y_3\Big{(}\int_{t_i}^{t_{i+1}}N\ind{X_1^N(s)\geq 1}\lambda_2(X_1^N(s),X_2^N(s))ds\Big{)} \\
 & \leq &  Y_3\Big{(}\int_0^T N\ind{X_1^N(s)\geq 1}\lambda_2(X_1^N(s),X_2^N(s))ds\Big{)} \\
\end{eqnarray*}
so by a similar argument as above, we also get the tightness of the BV norm
\begin{equation*}
 \prob{\|X_2^N\|_{[0,T]} \geq K } \to 0
\end{equation*}
as $ K\to 0$, independently in $N$. Then $X_2^N$ is tight in $L^p([0,T])$, for any $1\leq p < \infty$ (\cite{Giusti1984}).

\subsection{Identification of the limit}\label{ssec:proof_disc3}

We choose an adherence value $(0,X_2(t))$ of the sequence $(X_1^N(t),X_2^N(t))$ in $L^1([0,T])\times L^p([0,T])$. Then a subsequence (again denoted by) $(X_1^N(t),X_2^N(t))$ converge to $(0,X_2(t))$, almost surely and for almost $t \in [0,T]$. We are looking for test-functions such that
\begin{multline*}
 f(X_1^N(t),X_2^N(t))-f(X_1^N(0),X_2^N(0)-\int_0^t\mathbb{B}_N f(0,X_2^N(s))1_{X_1^N(s)=0}ds\\
 -\int_0^t\mathbb{B}_N f(X_1^N(s),X_2^N(s))1_{X_1^N(s)\geq1}ds
\end{multline*}
is a martingale and $\mathbb{B}_N f(X_1^N(s),X_2^N(s))$ is bounded independently of $N$ when $X_1\geq 1$. The following choice is inspired by \cite{Radulescu}. We introduce the stochastic process $Y_t^{x,y}$, starting at $y$ and whose generator is
\begin{equation*}
A^xg(y)= \lambda_2(x,y)\Big{[}g(y+1)-g(y)\Big{]},
\end{equation*}
for any $x\geq 1$. and we introduce the semigroup $P_t^x$ defined on bounded function, for any $x\geq 1$, by
\begin{equation}\label{eq:semigroupsize}
 P_t^x g(y)=\E\Big{[}g(Y_t^{x,y})e^{-\int_0^t\gamma_1(x,Y_s^{x,y})ds}\Big{]}.
\end{equation}
Then the semigroup $P_t^x$ satisfies the equation
\begin{equation*}
 \frac{dP_t^xg(y)}{dt}=A^xP_t^xg(y)-\gamma_1(x,y)P_t^xg(y).
\end{equation*}
Now for any bounded function $g$, define recursively
 \begin{align*}
 f(0,y)&=g(y),\\
 f(x,y)&=\int_0^{\infty}P_t^x(\gamma_1(x,.)f(x-1,.))(y)dt.
 \end{align*}
Such a test function is well defined by the assumption on $\gamma_1$. We then verify that
 \begin{align*}
  \B_Nf(0,y)&=\lambda_1(0,y)\Big{(}\int_0^{\infty}P_t^1(\gamma_1(1,.)g(.))(y)dt-g(y)\Big{)}+\gamma_2(0,y)\Big{[}g(y-1)-g(y)\Big{]}, \\
 \B_N f(x,y)&=\lambda_1(x,y)\Big{[}f(x+1,y)-f(x,y)\Big{]}+\gamma_2(x,y)\Big{[}f(x,y-1)-f(x,y)\Big{]}.
 \end{align*}
Indeed, for any $x\geq 1$,
\begin{eqnarray*}
& &A^xf(x,y)-\gamma_1(x,y)f(x,y) \\
&=&\int_0^{\infty}A^xP_t^x(\gamma_1(x,.)f(x-1,.))(y)-\gamma_1(x,y)P_t^x(\gamma_1(x,.)f(x-1,.))(y)dt,\\
&=&\int_0^{\infty}\frac{d}{dt}P_t^x(\gamma_1(x,.)f(x-1,.))(y)dt,\\
&=&\lim_{t \to \infty}P_t^x(\gamma_1(x,.)f(x,.))(y)-\gamma_1(x,y)f(x-1,y),\\
&=&-\gamma_1(x,y)f(x-1,y).
\end{eqnarray*}
Then
\begin{equation*}
 \lambda_2(x,y)\Big{[}f(x,y+1)-f(x,y)\Big{]}+\gamma_1(x,y)\Big{[}f(x-1,y)-f(x,y)\Big{]}=0.
\end{equation*}
Hence $\B_N f(x,y)$ is independent of N, and, taking the limit $N \to \infty$ in 
$$f(X_1^N(t),X_2^N(t))-f(X_1^N(0),X_2^N(0))-\int_0^t\mathbb{B}_N f(X_1^N(s),X_2^N(s))ds,$$
we deduce
\begin{equation*}
 g(X_2(t))-g(X_2(0))-\int_0^t\B_{\infty}g(X_2)
\end{equation*}
is a martingale where
\begin{equation*}
 \B_{\infty}g(y)=\lambda_1(0,y)\Big{(}\int_0^{\infty}P_t^1(\gamma_1(1,.)g(.))(y)dt-g(y)\Big{)}+\gamma_2(0,y)\Big{[}g(y-1)-g(y)\Big{]}. 
\end{equation*}

\paragraph*{Uniqueness}

Due to assumption on $k_1$ and $k_2$, the limiting generator defines a pure-jump Markov process in $\N$ which is not explosive. Uniqueness of the martingale then follows classically.
%

\section{Continuous-state bursting model}

The model we consider now is a continuous state-space model that explicitly assume the production of several molecules \textit{instantaneously}, through a jump Markov process, in agreement with experimental observations (\cite{Golding2005,Raj2006}). In line with experimental observations, it is standard to assume a Markovian hypothesis (an exponential waiting time between production jumps) and that the jump sizes are exponentially distributed (geometrically in the discrete case) as well (\cite{Suter2011}). The intensity of the jumps can be a linearly bounded function,  to allow for self-regulation. 

For simplicity, we will only consider the standard model of gene expression, that is with linear degradation rates and the production rate of the second variable is linear with respect to the first variable (as in example \ref{ex:1}). Note that more general rates as in the previous section could be considered as well. Here, we ask the question of what is the correct scaling so that the bursting production term is transmitted from the first variable to the second one, when the first variable is eliminated through an adiabatic limit. The propagation of bursting property in a gene network is an important question in molecular biology \cite{kaern05}.

This section is structured as follows. We first present the model in the rest of this paragraph, then state the results in subsection \ref{ssec:result_cont}, and divide the proofs in the remaining three subsection \ref{ssec:proofcont1},\ref{ssec:proofcont2},\ref{ssec:proofcont3}.

Let $x_1$ and $x_2$ denote the concentrations of mRNA and protein respectively. A simple model of single gene expression with bursting transcription is given by
\begin{eqnarray}
\label{eq:gene1}
\dfrac{dx_1}{dt}&=&-g_1 x_1 + \mathring{N}(h, \lambda_1(x_2)),\\
\label{eq:gene2}
\dfrac{dx_2}{dt}&=&-g_2 x_2+k_2 x_1.
\end{eqnarray}
Here $g_1$ and $g_2$ are the degradation rates for the mRNA and protein respectively, $k_2$ is the mRNA translation rate, and $\mathring{N}(h, \lambda_1(x_2))$ describes the transcription that is assumed to be a compound Poisson \textit{white noise} occurring at a rate $\lambda_1(x_2)$ with a non-negative jump size $\Delta x_1$ distributed with density $h$.

The equations (\ref{eq:gene1})-(\ref{eq:gene2}) are a short hand notation for
\begin{eqnarray}
\label{eq:gene1_sde}
x_1(t)&=&x_1(0)-\int_0^t g_1 x_1(s^-)ds \\
\nonumber & & +\int_0^t\int_{0}^{\infty}\int_{0}^{\infty}1_{\{r\le \lambda_1(x_2(s^-))\}} z N(ds,dz,dr), \\
\label{eq:gene2_sde}
x_2(t)&=&x_2(0)-\int_0^t g_2 x_2(s^-)ds+\int_0^t k_2 x_1(s^-)ds.
\end{eqnarray}
where $X_{s-}=\lim_{t\to s^-}X(t)$, and $N(ds,dz,dr)$ is a Poisson random measure on $(0,\infty)\times
[0,\infty)^2$ with intensity $ds h(z) dz dr$, where $s$ denotes the times of the jumps, $r$ is the state-dependency in an acceptance/rejection fashion, and $z$ the jump size.  Note that $(x_1(t))$ is a stochastic process with almost surely finite variation on any bounded interval $(0,T)$, so that the last integral is well defined as a Stieltjes-integral.

\begin{hyp}\label{hyp:funcadiab}
 The following discussion is valid for general rate functions $\lambda_1$ and density functions $h(\cdot)$ that satisfy
\begin{itemize}
\item $\lambda_1 \in C^1$, $\lambda_1$ is globally lipschitz and linearly bounded with 
\begin{equation*}
 0\leq \lambda_1(x)\leq c +K x.
\end{equation*}
\item $h\in C^0$ and $\int_{0}^\infty xh(x) dx < \infty$.
\end{itemize}
\end{hyp}
For such a general density function $h$, we denote the average burst size by
\begin{equation}
\label{eq:b}
b = \int_0^\infty x h(x)d x.
\end{equation}

\begin{remark}
Hill functions are often used to model gene self-regulation, so that $\lambda_1$ is given by
\begin{equation*}
 \lambda_1(x_2)=\frac{1+x_2^{\alpha}}{L+D x_2^{\alpha}}
\end{equation*}
where  $L$, $D$ are positive parameters and $\alpha$ is a positive integer (see \cite{Mackey2011} for more details). An exponential distribution of the bursting transcription is often used in modeling gene expression, in accordance with experimental findings (\cite{Xie2008}), so that the density function $h$ is given by
\begin{equation*}
h(x) = \dfrac{1}{b}e^{-x/b},
\end{equation*}
with $b$ the average burst size.
\end{remark}

If $\lambda_1(x_2) \equiv k_1$ is independent of the state $x_2$, the average transcription rate is
$bk_1$, and the asymptotic average mRNA and protein concentrations are
\begin{eqnarray*}
x_1^{\mathrm{eq}}:=\esp{x_1(t\to \infty)} &=& \frac{bk_1}{g_1}.\\
x_2^{\mathrm{eq}}:=\esp{x_2(t\to \infty)} &=& \frac{k_2}{g_2}x_1^{\mathrm{eq}} = \frac{bk_1k_2}{g_1g_2}.
\end{eqnarray*}

\subsection{Statement of the results}\label{ssec:result_cont}

Although Equations (\ref{eq:gene1})-(\ref{eq:gene2}) are simple, they are not analytically solvable. Hence, for pratical use to interpret experimental data, and to avoid numerical simulations, one looks for a reduced, analytically solvable, one-dimensional equation. In the following discussion, we consider the situation when mRNA degradation is a fast process, {\it i.e.} $g_1$ is ``large enough``, but the average equilibrium protein concentration $x_2^{\mathrm{eq}}$ remains unchanged. In most organisms and for most genes, mRNA has a smaller lifetime than protein (\cite{Schwanhausser2011}). In what follows, we denote by $g_1^n$, $k_2^n$ sequences of parameters, $\lambda_1^n$ sequence of functions and $h^n$ sequence of density function that will replace $g1$, $k_2$, $\lambda_1$, $h$ in \eqref{eq:gene1_sde}-\eqref{eq:gene2_sde} and satisfy hypothesis \ref{hyp:funcadiab}. We then denote $(x_1^n,x_2^n)$ its associated solution. We will always assume one of the following three scaling relations:
\begin{enumerate}
\item[(S1)] Frequent production rate of mRNA, namely $g_{1}^n=ng_1$, $\lambda_1^n=n\lambda_1$, and $k_2^n=k_2$ $h^n=h$ are independent of $n$;
\item[(S2)] Large burst of mRNA, namely $g_{1}^n=ng_1$, $h^n(z)=\frac{1}{n}h(\frac{z}{n})$ and $\lambda_1^n=\lambda_1$,$k_2^n=k_2$ remain unchanged;
\item[(S3)] Large production rate of protein, namely $g_{1}^n=ng_1$, $k_2^n=nk_2$, and $\lambda_1^n=\lambda_1$ $h^n=h$ are independent of $n$;
\end{enumerate}
These three different scaling are then associated with different behaviors of the biological systems given by $(x_1,x_2)$. As different genes may have different kinetics, each one of the possible scaling are of importance (\cite{Suter2011,Schwanhausser2011}). 

In this section we determine an effective reduced equation for equation \eqref{eq:gene2} for each of the three scaling conditions (S1)-(S3). In particular, we show that under assumption (S1), equation \eqref{eq:gene2} can be approximated by the deterministic ordinary differential equation
\begin{equation}
\label{eq:odey0}
\dfrac{dx_2}{dt}=-g_2 x_2+\lambda_2(x_2)
\end{equation}
where
\begin{equation*}
\lambda_2(x_2) = bk_2\lambda_1(x_2)/g_1.
\end{equation*}
We further show that under the scaling relations (S2) or (S3), equation \eqref{eq:gene2} can be reduced to the stochastic differential equation
\begin{equation}
\label{eq:gene2_reduce}
\dfrac{dx_2}{dt}=-g_2 x_2 + \mathring{N}(\bar{h},\lambda_1(x_2)).
\end{equation}
where $\bar{h}$ is a suitable density function in the jump size $\Delta x_2$ (to be detailed below).

We first explain, using some heuristic arguments, the differences between the three scaling relations and the associated results. When $n \to \infty$, $g_1^n \to \infty$ and applying  a standard quasi-equilibrium assumption we have
$$\frac{dx_1^n}{dt}\approx 0,$$
which yields
\begin{equation*}
 x_1^n(t) \approx \frac{1}{g_1^n}\mathring{N}(h^n(.),\lambda_1^n (x_2^n)) = \mathring{N}(g_1^n h^n(g_1^n\cdot), \lambda_1^n(x_2^n)),
\end{equation*}
and therefore the second equation (\ref{eq:gene2}) becomes
\begin{eqnarray*}
 \frac{dx_2^n}{dt} &\approx& -g_2 x_2^n+ \frac{k_2^n}{g_1^n}\mathring{N}(h^n(.),\lambda_1^n(x_2^n)) , \\
&\approx& -\gamma_2 x_2^n +\mathring{N}\left (\frac{g_1^n}{k_2^n}h^n(\frac{g_1^n}{k_2^n} ),\lambda_1^n(x_2^n)\right).
\end{eqnarray*}
Hence in \eqref{eq:gene2_reduce}, $\bar{h}(x_2) = (k_2/g_1)^{-1}h((k_2/g_1)^{-1}x_2)$ under the scaling (S2) and (S3). Furthermore, we note that the scaling (S2) also implies $n h^n(n \cdot) = h(\cdot)$, while in (S1), $n h^n(n \cdot) =n h(n\cdot)$
so that the jumps become more frequent and smaller.

We denote $(D[0,\infty),S)$ the cadlad function space of function defined on $[0,\infty)$ at values in $\R^+$ with the usual Skorohod topology (\cite{Jacod1987}). Similarly $(D[0,T],J)$ is the cadlag funtion space on $[0,T]$, with the Jakubowski topology (\cite{Jakubowski1997}). Also, $L^p[0,T)$ the space of $L^p$ integrable function on $[0,T)$, with $T>0$, which we endowed with total variation norm (\cite{Giusti1984}), and $M(0,\infty)$ is the space of real measurable function on $[0,\infty)$ with the metric (\cite{Kurtz1991})
\begin{equation*}
 d(x,y)=\int_0^{\infty}e^{-t}\max\{1,\mid x(t)-y(t) \mid \}dt. 
\end{equation*}
Our main results can be stated as follows
\begin{theorem}
 \label{th:1}
 Consider the equations \eqref{eq:gene1}-\eqref{eq:gene2} and assume Hypothesis \ref{hyp:funcadiab}.  If the scaling (S1) is satisfied, {\it i.e.}, $k_1^n  = nk_1$, and if $x_2^n(0)\to x_2^0$, then when $n\to\infty$,
 \begin{enumerate}
 \item The stochastic process $x_1^n(t)$ does not converge in any functional sense;
 \item The stochastic process $x_2^n(t)$ converges in law in $(D[0,\infty),S)$ towards the deterministic solution of the ordinary differential equation
\begin{equation}
\label{eq:ode}
\dfrac{dx_2}{dt}=-g_2 x_2+\lambda_2(x_2),\quad x_2(0)=x_2^0,
\end{equation}
where
\begin{equation*}
\lambda_2(x_2) = bk_2\lambda_1(x_2)/g_1.
\end{equation*}
\end{enumerate}
\end{theorem}

\begin{theorem}
 \label{th:2}
Consider the equations \eqref{eq:gene1}-\eqref{eq:gene2} and assume Hypothesis \ref{hyp:funcadiab}. If the scaling (S2) is satisfied, {\it i.e.}, $h^n(z)=\frac{1}{n}h(\frac{z}{n})$, and if $x_2^n(0)\to x_2^0$, then when $n\to\infty$,
\begin{enumerate}
\item The stochastic process $\frac{x_1^n(t)}{n}$ converges in law in $L^p$, $1\leq p < \infty$ and in $(D[0,T],J)$ to the (deterministic) fixed value $0$;
\item The stochastic process $x_2^n(t)$ converges in law in $L^p$, $1\leq p < \infty$ and in $(D[0,T],J)$ to the stochastic process defined by the solution of the stochastic differential equation
\begin{equation}
\label{eq:sde2}
\dfrac{dx_2}{dt}=-g_2 x_2+\mathring{N}(\bar{h},\lambda_1), \quad x_2(0) = x_2^0 \geq 0,
\end{equation}
where $\bar{h}(x_2) = (k_2/g_1)^{-1}h((k_2/g_1)^{-1}x_2)$.
\end{enumerate}
\end{theorem}

\begin{theorem}
 \label{th:3}
Consider the equations \eqref{eq:gene1}-\eqref{eq:gene2} and assume Hypothesis \ref{hyp:funcadiab}. If the scaling (S3) is satisfied, {\it i.e.}, $k_2^n  = nk_2$, and if $x_2^n(0)\to x_2^0$, then when $n\to\infty$,
\begin{enumerate}
\item The stochastic process $x_1^n(t)$ converges in law in $L^p$, $1\leq p < \infty$ and in $(D[0,T],J)$ to the (deterministic) fixed value $0$;
\item The stochastic process $x_2^n(t)$ converges in law in $L^p$, $1\leq p < \infty$ and in $(D[0,T],J)$ to the stochastic process determined by the solution of the stochastic differential equation
\begin{equation*}
\dfrac{dx_2}{dt}=-\gamma_2 x_2+\mathring{N}(\bar{h},\varphi), \quad x_2(0) = x_2^0 \geq 0,
\end{equation*}
where $\bar{h}(x_2) = (k_2/g_1)^{-1} h ((k_2/g_1)^{-1}x_2)$.
\end{enumerate}
\end{theorem}

\begin{remark}
Note that scalings (S2) and (S3) give similar results for the equation governing the protein variable $x_2(t)$ but very different results for the asymptotic stochastic process related to the mRNA. In particular, in Theorem \ref{th:2}, very large bursts of mRNA are transmitted to the protein, where in Theorem \ref{th:3}, very rarely is  mRNA present but when present it is efficiently synthesized into a burst of protein.
\end{remark}
In the rest of this paper, we provide proofs of the results mentioned above, using martingale techniques. In a companion paper \cite{Yvinec2012}, we use partial differential techniques to prove similar results (see also \cite{Haseltine2005,Zeron2010,Qian2011}).


The proofs of the three theorems above are divided in three steps. In section \ref{ssec:proofcont1} we first recall generator properties and derive moment estimates associated to \eqref{eq:gene1}-\eqref{eq:gene2}. In section \ref{ssec:proofcont2} we show the tightness result for all three theorems. We then identify the limit using a martingale approach in section \ref{ssec:proofcont3}.

%
\subsection{General properties and moment estimates}\label{ssec:proofcont1}

We first summarize the important background results on the stochastic processes used in the next.

\paragraph{One dimensional equation}
For the one-dimensional stochastic differential equation (\ref{eq:gene2_reduce}) perturbed by a compound Poisson white noise, of intensity $\lambda(.)$ and jump size distribution $\overline{h}(.)$, the extended generator of the stochastic process $(x_2(t))_{t\geq 0}$ is  \cite[Theorem 5.5]{Davis1984}, for any $f \in \mathcal{D}(\mathcal{A_1})$,
\begin{equation*}
 \mathcal{A}_1f(x) = -g_2 x\dfrac{df}{dx}+\lambda(x)\Big{(}\int_x^{\infty}\overline{h}(z-x)f(z)dz-f(x)\Big{)} 
\end{equation*}
\begin{eqnarray*}
\mathcal{D}(\mathcal{A}_1) &=&\{ f \in \mathcal{M}(0,\infty):\,\,t\mapsto f(xe^{-\gamma_2 t}) \text{ is absolutely }  \\
&& \text{ continuous for }t \in \mathcal{R}^+  \text{ and } \\
&& \ \mathbb{E}\sum_{T_i\leq t}|f(x_2(T_i))-f(x_2(T_i^-))|< \infty\text{ for all }t\geq0 \}
\end{eqnarray*}
where $\mathcal{M}(0,\infty)$ denotes a Borel-measurable function of $(0,\infty)$ and the times $T_i$ are the instants of the jump of $x_2$.
It is an extended domain containing all functions that are sufficiently smooth along the deterministic trajectories between the jumps, and with a bounded total variation induced by the jumps.

For any $f\in \mathcal{D}(\mathcal{A}_1)$, we have
\begin{equation*}
\frac{d\ }{dt}\mathbb{E}f(x_2(t))=\mathbb{E}\mathcal{A}_1(f(x_2(t))).
\end{equation*}

\paragraph{Two dimensional equation}  Consideration the two-dimensional stochastic differential equation \eqref{eq:gene1}-\eqref{eq:gene2} perturbed by a compound Poisson white noise, of intensity $\lambda_1(x_2)$ and jump size distribution $h$  follows along similar lines. Its infinitesimal generator and extended domain are
\begin{eqnarray}
\mathcal{A}_2g(x_1,x_2) &=& -g_1 x_1\dfrac{\partial g}{\partial x_1}+(k_2 x_1-g_2 x_2)\dfrac{\partial g}{\partial x_2} \nonumber \\
\label{generator2d} &&{} + \lambda_1(x_2)\Bigg{(}\int_{x_1}^{\infty}h(z-x_1)g(z,x_2)dz-g(x_1,x_2)\Bigg{)}, 
\end{eqnarray}
\begin{eqnarray}
\label{generator2d_b}
\mathcal{D}(\mathcal{A}_2) &=&\{ g \in \mathcal{M}((0,\infty)^2):\,\,t\mapsto g(\phi_t(x_1,x_2)) \text{ is absolutely }  \\
&& \text{ continuous for } t \in \mathcal{R}^+  \text{ and }\nonumber\\
&&\ \mathbb{E}\sum_{T_i\leq t}|g(x_1(T_i),x_2(T_i))-g(x_1(T_i^-),x_2(T_i^-))|< \infty\text{ for all }t\geq0 \} \nonumber
\end{eqnarray}
where $\phi_t$ is the deterministic flow given by the deterministic part of equations \eqref{eq:gene1}-\eqref{eq:gene2}, namely
\begin{eqnarray*}
\dfrac{dx_1}{dt}&=&-g_1 x_1,\\
\dfrac{dx_2}{dt}&=&-g_2 x_2+k_2 x_1.
\end{eqnarray*}

For any $f\in \mathcal{D}(\mathcal{A}_2)$, we have
\begin{equation}
\label{eq:weak}
\frac{d\ }{dt}\mathbb{E}f(x_1(t),x_2(t))=\mathbb{E}\mathcal{A}_2(f(x_1(t),x_2(t))).
\end{equation}



Using the stochastic differential equations \eqref{eq:gene1_sde}-\eqref{eq:gene2_sde}, we can deduce moment estimates, needed to be able to use unbounded test function (namely $f(x_1,x_2)=x_1$ and $f(x_1,x_2)=x_2$) in the martingale formulation. By taking the mean into \eqref{eq:gene1_sde}-\eqref{eq:gene2_sde}, neglecting negatives values and using hypothesis \ref{hyp:funcadiab},
\begin{eqnarray*}
 0&\leq& \esp{x_1(t)} \leq \esp{\int_0^t b\lambda_1(x_2(s))ds} \leq \int_0^t b(c +K\esp{x_2(s)})ds  \\
 0&\leq& \esp{x_2(t)} \leq \int_0^tk_2 \esp{x_1(s)}ds
\end{eqnarray*}
where we note $b=\esp{h}=\int_0^{\infty}zh(z)dz$. By Gronwall inequalities, there exist a constant $C$ such that
\begin{equation}\label{eq:estimMean}
\begin{aligned}
 \esp{\sup_{t\in [0,T]}x_1(t)}&\leq C(\esp{x_1(0)}+e^{CT})\\
\esp{\sup_{t\in [0,T]}x_2(t)}&\leq C(\esp{x_2(0)}+e^{CT})\\
\end{aligned}
\end{equation}
Then we claim that $f(x_1,x_2)=x_1$ is in the domain of the generator $\mathcal{A}_2$. We only have to verify (see Eq~\eqref{generator2d_b})
\begin{equation*}
  \mathbb{E}\sum_{T_i\leq t}|x_1(T_i)-x_1(T_i^-)|< \infty\text{ for all }t\geq0.
\end{equation*}
By equation \eqref{eq:gene1_sde}
\begin{eqnarray*}
\mathbb{E}\sum_{T_i\leq t}|x_1(T_i)-x_1(T_i^-)|&=& \mathbb{E}\int_0^t\int_{0}^{\infty}\int_{0}^{\infty}
1_{\{r\le \lambda_1(x_2(s^-))\}} z N(ds,dz,dr), \\
&\leq& b \esp{\int_0^t c +Kx_2(s) ds}.
\end{eqnarray*}
which is finite according to the previous estimates.

\subsection{Tightness}\label{ssec:proofcont2}

\paragraph*{S1}
We first show the tightness property for the scaling (S1) corresponding to theorem \ref{th:1}. In such case $x_1^n$ does no converge in any functional sense because it fluctuates very fast, as more and more jumps appears of size that stay of order $1$ (given by $h$). However, $\esp{x_1^n(t)}$ remains bounded, $\frac{x_1^n}{n}$ goes to $0$,
and by eq. \eqref{eq:gene2_sde},
\begin{equation*}
 \mid x_2^n(t) \mid \leq \mid x_2^n(0) \mid +\int_0^tk_2 \mid x_1^n(s)\mid ds.
\end{equation*}
For any $n$, let $N_n$ be the compound Poisson process associated to \eqref{eq:gene1_sde}, with $\{T_{n,i}\}_{i=1}^{\infty}$ the jump times which occur at a rate $n\lambda_1(x_2(s)^n)$, and $\{Z_{n,i}\}_{i=1}^{\infty}$ the jump sizes that are iid random variables with density $h$ (with the convention $T_{n,0}=0$ and $Z_{n,0}=X_0$),
\begin{equation*}
 N_n(t)=\sum_{T_{n,i}\leq t} Z_{n,i}.
\end{equation*}
 Then
\begin{equation*}
 x_1^n(t) = \sum_{T_{n,i}\leq t} Z_{n,i}e^{-ng_1 (t-T_{n,i})}\, \ind{t\geq T_{n,i}}.
\end{equation*}
By integration,
\begin{equation*}
 \int_0^t x_1^n(s)ds = \sum_{T_{n,i}\leq t} Z_{n,i}\frac{1}{ng_1}(1-e^{-g_1 (t-T_{n,i})})\, \ind{t\geq T_{n,i}}.
\end{equation*}
Then,
\begin{equation*}
 x_2^n(t) \leq x_2^n(0) +\int_0^t k_2 x_1^n(s)ds \leq x_2^n(0)+\frac{k_2}{ng_1}\sum_{T_{n,i}\leq t} Z_{n,i}.
\end{equation*}
Finally we deduce, by definition of the compound Poisson process,
\begin{equation*}
 x_2^n(t) \leq x_2^n(0)+\frac{k_2}{ng_1}N_n(t).
\end{equation*}
Now, by a time change, there exists a process $Y$ such that 
\begin{equation*}
N_n(t)=Y\big{(}\int_0^t n\lambda_1(x_2^n(s))ds \big{)}, 
\end{equation*}
where $Y$ is unit rate compound Poisson process of jump size iid (with density $h$). By the law of large number, $\frac{1}{n}Y(nt)$ is almost surely convergent, and hence almost surely bounded.
We deduce then there exists a random variable $C$ such that
\begin{equation*}
 x_2^n(t) \leq  x_2^n(0)+ \frac{k_2}{g_1}C\int_0^t \lambda_1(x_2^n(s))ds .
\end{equation*}
By Gronwall lemma and Markov inequality
\begin{equation*}
 \prob{\sup_{t \in [0,T]} x_2^n(t) \geq M } \to 0,
\end{equation*}
as $M\to \infty$ and uniformly in $n$. Similarly, for any $t_1,t_2 \in [0,T]$, 
\begin{equation*}
 \mid x_2^n(t_2)-x_2^n(t_1) \mid \leq \frac{k_2}{ng_1} \mid N_n(t_2)-N_n(t_1) \mid .
\end{equation*}
Again, $\displaystyle{ N_n(t_2)-N_n(t_1)=Y\big{(}\int_{t_1}^{t_2} n\lambda_1(x_2^n(s))ds \big{)} }$ and, still by the law of large number
\begin{equation*}
 \mid x_2^n(t_2)-x_2^n(t_1) \mid \leq \frac{k_2}{g_1}C\int_{t_1}^{t_2} \lambda_1(x_2^n(s))ds,
\end{equation*}
so that , for any $\varepsilon>0$
\begin{equation*}
 \lim_{\theta \to 0}\limsup_n \sup_{S_1\leq S_2 \leq S_1+\theta} \prob{\mid x_2^n(S_2)-x_2^n(S_1)\mid \geq \varepsilon}=0,
\end{equation*}
where the supremum is over stopping times bounded by $T$. Then by Aldous' tightness criteria (\cite[thm 4.5 p 356]{Jacod1987}), $x_2^n$ is tight in $(D[0,\infty),S)$.

\paragraph*{S3}
Now we show the tightness property for the scaling (S3) corresponding to theorem \ref{th:3}, with $k_2^n=nk_2$. In such case $x_1^n$ converges to $0$ in $L^1$, and we get a control over $n\int_0^t x_1^n(s)ds$. 
Indeed using $g(x_1,x_2)=x_1$ in \eqref{generator2d}, we get
\begin{eqnarray*}
 x_1^n(t)-x_1^n(0)-\int_0^t (-ng_1 x_1^n(s)+ b\lambda_1(x_2(s)^n) ds),
\end{eqnarray*}
is a martingale so that due to Hypothesis \ref{hyp:funcadiab}, 
\begin{equation*}
g_1\esp{ n\int_0^t x_1^n(s)ds} \leq \esp{x_1^n(0)} + b(ct +K\int_0^t \esp{x_2^n(s)}ds )
\end{equation*}
By eq. \eqref{eq:gene2_sde},
\begin{eqnarray*}
  x_2^n(t)  &\leq&  \esp{x_2^n(0)}+k_2 n\int_0^t x_1^n(s)ds.
\end{eqnarray*}
then 
\begin{eqnarray*}
  \sup_{t\in [0,T]} x_2^n(t)  & \leq &  \esp{x_2^n(0)}+k_2 n\int_0^T x_1^n(s)ds.
\end{eqnarray*}
Reporting into the estimates for $x_1^n$ yelds 
\begin{eqnarray*}
g_1\esp{ n\int_0^t x_1^n(s)ds} &\leq &\esp{x_1^n(0)} + b (ct +K(\esp{x_2^n(0)}+tk_2 n\int_0^t \esp{x_1^n(s)}ds )),\\
&\leq& C_T^1 +C_T^2 \esp{n\int_0^t x_1^n(s)ds},
\end{eqnarray*}
for two constants $C_T^1$, $C_T^2$ that depends solely on $T$. Then $\esp{ n\int_0^t x_1^n(s)ds}$ is bounded uniformly in $n$ so that $x_1^n$ converges to $0$ in $L^1$ and
 \begin{equation*}
 \prob{\sup_{t \in [0,T]} x_2^n(t) \geq M } \to 0
\end{equation*}
as $M\to \infty$ and uniformly in $n$. Now for any subdivision of $[0,T]$, $0=t_0 <t_1<\cdots <t_n=T$,
\begin{eqnarray*}
 \sum_{i=0}^{n-1} \mid x_2^n(t_{i+1})-x_2^n(t_i) \mid &\leq & \esp{x_2^n(0)}+k_2 n\int_0^t x_1^n(s)ds,
\end{eqnarray*}
so that we also get the tightness of the BV norm,
\begin{equation*}
 \prob{\|x_2^n\|_{[0,T]} \geq M } \to 0,
\end{equation*}
as $ M\to 0$, independently in $n$. Then $x_2^n$ is tight in $L^p([0,T])$, for any $1\leq p < \infty$ (\cite{Giusti1984}) and also, by a similar criteria, in  $(D[0,T],J)$ (\cite{Jakubowski1997}).
 
\paragraph*{S2}
Now we show the tightness property for the scaling (S2) corresponding to theorem \ref{th:2}, with $h^n=\frac{1}{n}h(\frac{1}{n})$. Remark that on such case, denoting $z^n=\frac{x_1^n}{n}$, the variables $(z^n,x_2^n)$ satisfies \eqref{eq:gene1_sde}-\eqref{eq:gene2_sde} with the (S3) scaling, so we already now that $x_2^n$ is tight in $L^p([0,T])$, for any $1\leq p < \infty$.

For $x_1^n$, formally, note that each jumps yelds a contribution for $\int x_1^n$ of $\frac{b}{g_1}$ so there's no hope for a convergence to $0$ in $L^1$. However, we still have
$$x_1^n(t) = \sum_{T_{n,i}\leq t} Z_{n,i}e^{-ng_1 (t-T_{n,i})}\, 1_{t\geq T_{n,i}}.$$
where $T_{n,i}$ appears with rate $\lambda_1(x_2^n(s))$, and $\{Z_{n,i}\}_{i=1}^{\infty}$ are iid random variables with density $h^n$. Then
\begin{eqnarray*}
 x_1^n(t) &\leq& \sum_{T_{n,i}\leq t} Z_{n,i}\Big{(}\ind{[T_{n,i},T_{n,i}+\frac{1}{\sqrt{n}}}+e^{-ng_1 \frac{1}{\sqrt{n}}}\, 1_{t\geq T_{n,i}}\Big{)}.
\end{eqnarray*}
But for $M>0$, by Markov inequality,
\begin{equation*}
 \prob{Z_{n,i}e^{-\sqrt{n}g_1}>M}\leq \frac{nb}{Me^{\sqrt{n}g_1}}\leq \varepsilon,
\end{equation*}
for any $\varepsilon$ and $n$ sufficiently large. Then, conditionning by the jump times,
\begin{eqnarray*}
 \int_0^t \prob{x_1^n(s)>M \mid T_{n,i}}\leq \sum_{T_{n,i}\leq t} \frac{1}{\sqrt{n}}\ind{t\geq T_{n,i}}+ \sum_{T_{n,i}\leq t}\varepsilon (t-T_{n,i})\ind{t\geq T_{n,i}}\leq \varepsilon.
\end{eqnarray*}
for $n$ large. Because $\int_0^t x_2^n(s)ds$ has been shown to be bounded independently of $n$, we can drop the conditionning, and $\int_0^t \prob{x_1^n(s)>M }$ is arbitrary small. We show also similarly that  
\begin{equation*}
\lim_{h \to 0}\sup_n \int_0^T \max (1,\mid x_1^n(t+h)-x_1^n(t) \mid )dt=0,
\end{equation*}
so that $x_1^n$ is tight in $M(0,\infty)$ (\cite[thm 4.1]{Kurtz1991}).

\subsection{Identification with the martingale problem}\label{ssec:proofcont3}

The three theorems below can be proved using martingale techniques, with similar spirit. For each scaling, the generator $\mathcal{A}_2^n$ can be decomposed into a fast component, or order $n$, and a slow component, of order $1$. In each case, one need to find  particular condition to ensure that the fast component vanishes. For the scaling $(S1)$, the fast component acts only in the first variable, so ergodicity of this component will ensure that it vanishes, as in averaging theorems \cite{Kurtz1992}. For the other two scaling, the fast component acts on both variables, and we will have to find the particular relation between both variable that ensures this component vanishes.

\paragraph{Proof of Theorem \ref{th:1}}

For any $B \in \B(\R_+)$, $t>0$, we define the occupation measure 
\begin{equation*}
V_1^n(B\times [0,t])=\int_0^t \ind{B}(x_1^n(s))ds,
\end{equation*}
and we identify $V_1^n$ as a stochastic process with value in the space of finite meaure on $\R^+$. Because $\esp{x_1^n(t)}$ remains bounded uniformly in $n$ on any $[0,T]$, it is stochastically bounded and $V_1$ then satisfies Aldous criteria of tightness. Now take a test function $f$ that depends only on $x_1$, so that
\begin{equation*}
 \mathcal{A}_2^nf(x_1)=nC_{x_2}f(x_1),
\end{equation*}
with
\begin{equation*}
 C_{x_2}f(x_1)=-g_1 x_1 f'(x_1) + \lambda_1(x_2)\Bigg{(}\int_{x_1}^{\infty}h(z-x_1)f(z)dz-f(x_1)\Bigg{)}.
\end{equation*}
Then
\begin{equation*}
 M_t^n=f(x_1^n(t))-f(x_1^n(0))-n\int_{\R_+}\int_0^t C_{x_2^n(s)}f(x_1)V_1^n(dx_1 \times ds)
\end{equation*}
is a martingale. Dividing by $n$, for any limiting point $(V_1,x_2)$, we must have, for any $f\in C_b(\R_+)$,
\begin{equation*}
 \esp{\int_{\R_+}\int_0^t C_{x_2(s)}f(x_1)V_1(dx_1 \times ds)}=0.
\end{equation*}
Because for any $x_2$, the generator $C_{x_2}$ is (exponentially) ergodic, $V_1$ is uniquely determined by the invariant measure associated to $C_{x_2}$. In particular, for any $t>0$
\begin{equation*}
 \int_{\R_+}\int_0^t x_1 V_1^n(dx_1 \times ds) \to \int_0^t \frac{b}{g_1}\lambda_1(x_2(s))ds.
\end{equation*}
Then for $f$ that depends only on $x_2$,
\begin{equation*}
 f(x_2^n(t))-f(x_2^n(0))-\int_{\R_+}\int_0^t (k_2 x_1-g_2 x_2^n(s))f'(x_2^n(s))V_1^n(dx_1 \times ds)
\end{equation*}
converges to
\begin{equation*}
 f(x_2(t))-f(x_2(0))-\int_0^t (\frac{bk_2}{g_1}\lambda_1(x_2(s))-g_2 x_2(s))f'(x_2(s))ds
\end{equation*}
Due to the assumption on $\lambda_1$, there exists a unique solution associated to the (deterministic) equation \ref{eq:odey0} so $x_2$ is uniquely determined.

\paragraph{Proof of Theorem \ref{th:3}}

We've already seen that $x_1^n$ converges to $0$ in $L^1([0,T])$ and $x_2^n$ is tight in $L^p([0,T])$. We then take a subsequence  $(x_1^n(t),x_2^n(t))$  that converges to $(0,x_2(t))$, almost surely and for almost $t \in [0,T]$. Then we consider the fast component of the generator $\mathcal{A}_2^n$, given in this case by
\begin{equation*}
 -g_1 x_1 \frac{\partial f}{\partial x_1}+k_2 x_1\frac{\partial f}{\partial x_2}.
\end{equation*}
This defines a transport equation. Starting at $(x_1,x_2)$ at time $0$, the asymptotic value of the flow associated to the transport equation is $(0,y)$ where
\begin{equation*}
 y=x_2+\int_0^{\infty}k_2 x_1(s)ds=x_2 +\int_0^{x_1}\frac{k_2 z}{g_1 z}dz=x_2+\frac{k_2}{g_1}x_1
\end{equation*}
We then consider
\begin{equation*}
 f(x_1,x_2)=g(x_2+\frac{k_2 }{g_1 }x_1),
\end{equation*}
that satisfies, for any $x_1,x_2$,
\begin{equation*}
 -g_1 x_1 \frac{\partial f}{\partial x_1}+k_2 x_1\frac{\partial f}{\partial x_2}=0.
\end{equation*}
Now taking the limit $n\to \infty$ into
\begin{equation*}
 f(x_1^n(t),x_2^n(t))-f(x_1^n(0),x_2^n(0))-\int_0^t \mathcal{A}_2^nf(x_1^n(s),x_2^n(s))ds,
\end{equation*}
yelds
\begin{multline*}
  g(x_2(t))-g(x_2(0))-\int_0^t -g_2 x_2g'(x_2(s)) \\
+\lambda_1(x_2(s))\Bigg{(}\int_{0}^{\infty}\bar{h}(z)g(x_2(s)+z)dz-g(x_2(s))\Bigg{)}ds,
\end{multline*}
where $\bar{h}(x_2) = (k_2/g_1)^{-1}h((k_2/g_1)^{-1}x_2)$. Hence the limiting process $x_2$ must satisfy the martingale problem associated with the generator
\begin{equation*}
 \mathcal{A}_{\infty}g(x) = -g_2 x\dfrac{dg}{dx}+\lambda_1(x)\Big{(}\int_x^{\infty}\bar{h}(z-x)f(z)dz-f(x)\Big{)}, 
\end{equation*}
for which uniqueness holds for bounded $k_1$ (see \cite[thm 2.5]{Radulescu}). A truncature argument allows then to conclude for linearly bounded $k_1$.

\paragraph{Proof of Theorem \ref{th:2}}

As noticed before, $(z^n,x_2^n)$ with $z^n(t)=\frac{x_1^n(t)}{n}$ satisfies the scaling (S3) so similar conclusion holds for $x_2^n$.

\bibliographystyle{imsart-nameyear}
\bibliography{FastSlowVar_Jan31_MCM}

\end{document}